\batchmode
\documentclass [11pt]{article}
\usepackage[latin1]{inputenc} 
\usepackage{amssymb}

\newtheorem{theorem}{Theorem}
\newtheorem{corollary}[theorem]{Corollary}

\newtheorem{remark}[theorem]{Remark}
\newtheorem{definition}[theorem]{Definition}
\newtheorem{example}[theorem]{Example}

\title{ Rotation Minimizing  vector fields and  frames\\ in Riemannian manifolds}

\author{Fernando Etayo\footnote{Departamento de Matem\'{a}ticas, Estad\'{\i}stica y
Computaci\'{o}n.
 Facultad de Ciencias.  Universidad de Cantabria.
 Avda. de los Castros, s/n, 39071 Santander, SPAIN.
 e-mail: etayof@unican.es}}
\date{}

\begin{document}
\maketitle

\begin{abstract}

We prove that a normal vector field along a curve in $\mathbb{R}^{3}$  is rotation minimizing (RM) if and only
if it is parallel  respect to the normal connection. This allows us to generalize all the results of RM vectors
and frames to curves immersed in Riemannian manifolds.

\end{abstract}



\section{Introduction}

In a celebrated paper \cite{B}, Bishop introduced what nowadays are called rotation minimizing vector fields
(RM, for short) over a curve in the Euclidean 3-space. The purpose of this note is to show how they can be
defined in Riemannian manifolds.

When one considers vector fields and moving frames along a curve, one can take into account two  general ideas:
defining an adapted frame whose first vector is the tangent vector to the curve and defining a frame whose first
vector is the position vector of the point of the curve. In both cases, one can consider Frenet frames and RM
frames. See \cite{G}. The first idea, frames containing the tangent vector of the curve as first vector, can be
generalized to the case of a curve immersed in a Riemannian manifold, while the second one has no sense in this
general framework.

We shall introduce such definition of a RM vector field  over a curve in a Riemannian manifold as a vector field
parallel respect to the normal connection. It will be shown that this definition is consistent with that of
Bishop for curves in the Euclidean space (section 4), and that remains invariable under isometries (section 5).
In section 2 we shall remember the basic definitions about RM vector fields and frames, and in section 3 about
curves in Riemannian manifolds.

\section{RM vector fields and frames of a curve in $\mathbb{R}^{3}$}

\begin{definition} A normal vector field $\vec{v}=\vec{v}(t)$ over a curve $\gamma =\gamma (t)$ in $\mathbb{R}^{3}$
is said
to be  \emph{relatively parallel} or \emph{rotation minimizing} (RM) if the derivative $\vec{v}'(t)$ is
proportional to $\gamma '(t)$.

\end{definition}

\begin{remark}

(1) In this case the ruled surface $f(t,\lambda ) = \gamma (t) +\lambda \vec{v}(t)$ is developable, because
$[\gamma '(t), v(t), v'(t)]=0$.

(2) If $v$ is a RM vector field, then $\parallel \vec{v}\parallel$ is constant. Let $\vec{t}$ denote the tangent
vector to $\gamma$. Then $\vec{v}' = \lambda \vec{t} \Rightarrow \vec{v}' \perp \vec{v} \Rightarrow\frac{d}{dt}
(\vec{v}\cdot \vec{v})=0 $.

\end{remark}

\begin{example} (1) Let $\gamma (t)=t$ be the line given by the x-axis in $\mathbb{R}^{3}$.
Then a normal vector field $\vec{v}(t)$ over $\gamma $ is RM respect to $\gamma$  iff it is constant.

(2) Let $\gamma (t)=(\cos t,\sin t,0)$ the unit circle in the horizontal plane. Let $(0,0,h)$ any point in the
vertical axis. Let us consider the vector $\vec{v}(t)$ joining $x(t)$ and $(0,0,h)$. Then $\vec{v}$ is RM. The
developable surface generated is the corresponding cone.

(3) Normal and binormal vector fields, $ \vec{n}$ and $\vec{b}$, of a Frenet moving frame are not RM vector
fields in general. Let $\gamma =\gamma (s)$ a curve parametrized respect to the arc-length. Then, the
Frenet-Serret formulas say that $\vec{n}' = -\kappa \vec{t} + \tau \vec{b}$, and $\vec{b}=-\tau \vec{n}$ where
$\kappa$ and $\tau$ denote the curvature and the torsion. For  a twisted (non plane) curve, $\tau \neq 0$,   the
above equations show that $\vec{n}$ and $\vec{b}$ are nor RM vector fields. This is the case, for instance, of
the helix $(a\cos t, a\sin t, bt)$. The normal vector is $(-\sin t, \cos t,0)$ and the surface generated by the
normal lines is the helicoid, which is not developable, because its Gauss curvature does not vanish. For a plane
curve, $ \vec{n}$ and $\vec{b}$ are RM vector fields.
\end{example}

\begin{definition} Let $\gamma =\gamma (t)$ in $\mathbb{R}^{3}$ be a curve. A RM \emph{frame}, \emph{parallel frame},
\emph{natural frame}, \emph{Bishop frame} or \emph{adapted frame} is a moving orthonormal frame $\{
\vec{t}(t),\vec{u}(t),\vec{v}(t)\}$ along $\gamma$, where $\vec{t}(t)$ is the tangent vector to $\gamma$ at the
point $\gamma (t)$ and $\vec{u},\vec{v}$ are RM vector fields.
\end{definition}

\begin{example} If $\vec{u}$ is a unitary RM vector field along $\gamma$, then $\{ \vec{t},\vec{u},\vec{t}\times
\vec{u}\}$is a RM frame along $\gamma$.

\end{example}

RM frames are a very useful tool in some aspects of Computer Aided Geometric Design. See \cite{F} as a basic
reference.

\section{Curves immersed in a Riemmanian manifold}

In this section we shall remember the main facts about curves immersed in a Riemmanian manifold.

As it is well known, if $M$ is a submanifold of a Riemannian manifold $(\overline{M},\overline{g})$, then the
Levi-Civita connection $\overline{\nabla }$ of $(\overline{M},\overline{g})$ induces a Levi-Civita connection in
$(M,g=\overline{g}|_{M})$ and a normal connection $D^{\perp }:\mathfrak{X}(M)\times \Gamma TM^{\perp }\to \Gamma
TM^{\perp }$, where $\mathfrak{X}(M)$ denotes the module of vector fields of the submanifold and $\Gamma
TM^{\perp }$ the module of sections of the normal bundle. We shall use the notation of \cite[§ VIII]{KN}.

The normal connection is defined as follows. First of all, consider a point $p\in M$. Then the tangent space at
the manifold $\overline{M}$ can be decomposed as a orthogonal direct sum $T_{p}(\overline{M}) =T_{p}M \oplus
T^{\perp }_{p}M$, where $T_{p}M$ (resp. $T^{\perp }_{p}M$) denotes the tangent (resp. the normal) space to the
submanifold. The set $TM=\bigcup_{p\in M}T_{p}M$ (resp. $T^{\perp }M=\bigcup_{p\in M}T^{\perp }_{p}M$) is a
manifold called the tangent bundle (res. the normal bundle) and it has a a structure of vector bundle over $M$,
given by the natural projection $\pi :TM\to M; \pi (u)=p$ if $u\in T_{p}M$ (resp. $\pi ^{\perp }:TM^{\perp }\to
M; \pi (v)=p$ if $v\in T^{\perp }_{p}M$). Let us denote by $\mathfrak{X}(M)$ (resp. $\Gamma TM^{\perp }$) the
module of vector fields over $M$, i.e.,  the module of sections of $\pi :TM\to M$ (resp. the module of sections
of  $\pi ^{\perp }:TM^{\perp }\to M)$.

Then for any $X\in \mathfrak{X}(M)$ and $v\in \Gamma TM^{\perp }$ one has the  decomposition

$$\overline{\nabla }_{X}v = -A_{v}X+ D^{\perp }_{X}v $$

\noindent where $-A_{v}X\in \mathfrak{X}(M)$ and $D^{\perp }_{X}v \in \Gamma TM^{\perp }$. The \emph{Weingarten
operator} $A$ is a $\mathfrak{F}(M)$-bilinear map and the normal connection $D^{\perp }$ is a  connection in the
normal  bundle $T^{\perp }M\to M$.

Moreover, if $v,w \in \Gamma TM^{\perp }$ are two normal vector fields, then

$$\overline{g}(D^{\perp }_{X}v,w)+  \overline{g}(v,D^{\perp }_{X}w) =X(\overline{g}(v,w))
$$

\noindent which shows that the normal connection $D^{\perp }$ is metric for the fibre metric in the normal
bundle $TM^{\perp }$.

A normal vector field $v$ is said to be \emph{parallel} respect to $X\in \mathfrak{X}(M)$ if $D^{\perp
}_{X}v=0$.

Let $\pi :E\to M$ be a vector bundle with a connection  $D$, compatible with a metric on $E$.  The notion of
being a section $\sigma :M\to E$ parallel respect to the connection $D$ is equivalent to that it is defined by
the parallel transport induced by the connection. See \cite{P} for details. The parallel transport defines an
isometry between any two different fibres of $\pi :E\to M$. Thus, the norm of a parallel section remains
constant and the angle between two parallel sections also remains constant.  In the present case of having a
submanifold $M$ of a Riemannian manifold $(\overline{M},\overline{g})$,
 the vector bundle is the normal  bundle, and the metric is the restriction of
 $\overline{g}$ to normal vectors.

Many results about curves un Riemannian manifolds have been obtained in the past. We would like to to point out
that generalizations of Frenet formula have been obtained in \cite{M}, in \cite{MR} for the case of spaces of
constant curvature, and in \cite{OE} for the case of the Minkowski space. Besides in \cite{CMM} some results
about the total curvature of a curve in a Riemannian manifold are also obtained.

\section{RM vector fields over a curve immersed in a Riemmanian manifold}

The notion  of RM vector field implies a notion of parallel transport. We shall show this carefully.

 \bigskip

Let us consider the case where $M$ is a curve $\gamma$ and $(M,g=\overline{g}|_{M})$ is $\mathbb{R}^{3}$ with
the standard product. Let us denote by $T_{\gamma (t_{0})}$ (resp. $T^{\perp }_{\gamma (t_{0})}$) the tangent
line (resp. the normal plane) to $\gamma $ in $\gamma (t_{0})$. Then we have:

\begin{theorem} A normal vector field $v$ over a curve $\gamma$ immersed in $\mathbb{R}^{3}$ is a RM vector field
iff ir is parallel
respect to the normal connection of $\gamma$.
\end{theorem}

\emph{Proof.} Let us denote as $(x^{1},x^{2},x^{3})$ the global coordinates in $\mathbb{R}^{3}$. The curve
$\gamma$ can be expressed as $\gamma (s)=(\gamma ^{1}(s), \gamma ^{2}(s), \gamma ^{3}(s))$, $s$ being the
arc-length parameter,  and the tangent vector is $\overline{t}= \gamma '(s) = \frac{\partial }{\partial
x^{i}}\frac{d \gamma ^{i}}{ds}$ (Einstein's convention is assumed).

Let $v$ be a normal vector field over $\gamma$, $v=\frac{\partial }{\partial x^{i}}v^{i}$. The condition of
being normal to the curve means that

$$
 \sum_{i=1}^{3} \hspace{2mm} v^{i}\hspace{2mm}\frac{d \gamma ^{i}}{ds}=0
$$

The condition of being $v$ a RM vector field means that $v'(t)$ is proportional to $\gamma '(t)$, i.e.

$$\label{RM} \frac{d v^{i}}{ds} = \lambda (s)\hspace{2mm} \frac{d \gamma
^{i}}{ds}, \hspace{1cm} \forall i=1,2,3
 $$

\noindent where $\lambda =\lambda (s)$ is a function.

\bigskip

Now, we shall check the value of $D^{\perp }_{\gamma '(s)}v$. We must prove that $D^{\perp }_{\gamma '(s)}v=0$
iff the above equation  is satisfied. Let $\overline{\nabla }$ be the Levi-Civita connection of
$\mathbb{R}^{3}$. The all of its Christoffel symbols vanish, and then one has:

$$\overline{\nabla }_{\overline{t}}\hspace{2mm}v =
 \overline{\nabla }_{\frac{\partial }{\partial x^{i}}\frac{d \gamma ^{i}}{ds}} \hspace{2mm} \left( \frac{\partial }
 {\partial x^{j}}v^{j} \right) =
\frac{d \gamma ^{i}}{ds}\overline{\nabla }_{\frac{\partial }{\partial x^{i}}} \hspace{2mm} \left( \frac{\partial
}{\partial x^{j}}v^{j} \right) =
 \frac{\partial }{\partial x^{j}} \frac{d \gamma ^{i}}{ds}\frac{\partial v^{j} }{\partial x^{i}}
$$

Applying the chain rule one has:

$$\frac{\partial v^{j} }{\partial x^{i}} = \frac{dv^{j} }{ds}\frac{ds}{dx^{i}}=
\frac{dv^{j} }{ds} \left( \frac{dx^{i}}{ds}    \right) ^{-1} = \frac{dv^{j} }{ds} \left( \frac{d(x^{i}\circ
\gamma ^{-1})}{ds}    \right) ^{-1} =
 \frac{dv^{j} }{ds} \left( \frac{d\gamma ^{i}}{ds}    \right) ^{-1}.$$

And then,

$$ \overline{\nabla }_{\overline{t}}\hspace{2mm}v =
\frac{\partial }{\partial x^{j}} \frac{d \gamma ^{i}}{ds}\frac{\partial v^{j} }{\partial x^{i}} =
 \frac{\partial }{\partial x^{j}} \frac{d \gamma ^{i}}{ds} \frac{dv^{j} }{ds} \left( \frac{d\gamma ^{i}}{ds}
 \right) ^{-1}  =  \frac{\partial }{\partial x^{j}} \frac{dv^{j} }{ds}.
 $$

 Finally, $v$ is parallel $\Leftrightarrow D^{\perp }_{\gamma '(s)}v=0
 \Leftrightarrow \overline{\nabla }_{\overline{t}}v $ is tangent to $\gamma \Leftrightarrow \frac{d v^{i}}{ds} =
  \lambda (s)\hspace{2mm} \frac{d \gamma
^{i}}{ds}, \forall i=1,2,3  \Leftrightarrow v$ is a RM vector field, thus finishing the proof. $\Box$

\bigskip

The above result is important because it allows to obtain the definition of a RM vector field over a curve
immersed in a Riemmannian manifold. Moreover, one easily can deduce the following properties of the above
Proposition.

\begin{corollary} With the above notation:

(1) Given a vector $v_{0}\in T^{\perp }_{\gamma (t_{0})}$ there exists a unique RM vector field $v$ over
$\gamma$ such that $v(t_{0}) =v_{0}$.

(2) If $v$ is a RM vector fields over $\gamma$ then the norm $\| v\|$ is constant.

(3) If $v$ and $w$ are RM vector fields over $\gamma$ then the angle between $v(t)$ and $w(t)$ is constant.
$\Box$
\end{corollary}

Thus, we can give the following

\begin{definition} Let $\gamma$ be a curve immersed in a Riemannian manifold $(\overline{M},\overline{g})$.

(1) A normal vector field $v$ over $\gamma$ is said to be a \emph{RM vector field} if it is parallel respect to
the normal connection of $\gamma$.

(2) A \emph{parallel frame}, \emph{natural frame}, RM \emph{frame} or \emph{adapted frame} is a moving
orthonormal frame $\{ \vec{t}(t),\vec{v}_{1}(t),\ldots ,\vec{v}_{n}(t)\}$ along $\gamma$, where $\vec{t}(t)$ is
the tangent vector to $\gamma$ at the point $\gamma (t)$ and $\vec{v}_{i}$ are RM vector fields, $\forall i\in
{1,\ldots ,n}$.
\end{definition}

If one defines an orthonormal frame $\{ \vec{t}(t_{0}),\vec{v}_{1}(t_{0}),\ldots ,\vec{v}_{n}(t_{0})\}$ at a
point $\gamma (t_{0})$ of a curve $\gamma$ then by parallel transport  it can be extended along $\gamma$.
Parallel transport is an isometry, this meaning that norms and angles are preserved.

\section{RM frames and transformations}

First of all, we prove that RM vector fields and frames are preserved by isometries.  Let $\mu:
(\overline{M},\overline{g})\to (\overline{M},\overline{g})$ be an isometry and let $\mu _{* p}
:T_{p}\overline{M}\to T_{p}\overline{M}$ its differential or tangent map. Then $\mu _{* p}$ is a linear isometry
respect to $\overline{g}_{p}$, i.e., $\overline{g}(\mu _{*}v,\mu _{*}w) = \overline{g}(v,w)$.

\begin{theorem}  Let $\gamma$ be a curve immersed in a Riemannian manifold $(\overline{M},\overline{g})$ and let
$\mu: (\overline{M},\overline{g})\to (\overline{M},\overline{g})$ be an isometry.

(1) If  $v$ is a RM vector field over $\gamma$, then $\mu _{*}(v)$ is a RM vector field over $\mu \circ \gamma$.

(2) If $\{ \vec{t},\vec{v}_{1},\ldots \vec{v}_{n}\}$ is a RM frame over $\gamma$, then $\{ \mu _{*}(\vec{t}),\mu
_{*}(\vec{v}_{1}),\ldots ,\mu _{*}(\vec{v}_{n})\}$ is a RM frame over $\mu \circ \gamma$.

\end{theorem}

\emph{Proof.} The following claims are well known:

\begin{enumerate}
\item  If $\vec{t}(t_{0})$ is the tangent vector of $\gamma$ at the point $\gamma (t_{0})$, then $\mu
_{*}(\vec{t}(t_{0}))$ is the tangent vector of $\mu \circ \gamma$ at the point $(\mu \circ \gamma)(t_{0})$.

\item If $\vec{v}\in T^{\perp }_{\gamma ({t_{0}})}$, then $\mu_{*}(\vec{v})\in T^{\perp }_{(\mu \circ \gamma )
({t_{0}})}$, because $\mu _{*}$ is an isometry.

\item $\mu _{*}(\overline{\nabla }_{X}Y)=\overline{\nabla }_{\mu _{*}X}\hspace{2mm}\mu _{*}Y$ (cfr., e.g
\cite[p. 161, vol.1]{KN}).

\end{enumerate}

Claims (1) and (2) show that $\mu _{*}$ maps tangent (resp. normal) vectors in tangent (resp. normal vectors).
\bigskip

Now, we can easily prove the theorem.

(1) Let $v$ be a RM vector field over $\gamma$. Then $D^{\perp }_{\vec{t}} v=0$, where $\vec{t}$ denotes the
tangent vector to $\gamma$. We must prove that $D^{\perp }_{\mu_{*}\vec{t}}\hspace{2mm} \mu_{*}v=0$.

We have the tangent and normal decomposition:

$$\overline{\nabla }_{\mu_{*}\vec{t}}\hspace{2mm} \mu_{*}v = -A_{\mu_{*}v}\hspace{2mm}\mu_{*}\vec{t} +
D^{\perp }_{\mu_{*}\vec{t}}\hspace{2mm} \mu_{*}v
$$

\noindent and, on the other hand,

$$\overline{\nabla }_{\mu_{*}\vec{t}}\hspace{2mm} \mu_{*}v =\mu_{*}(\overline{\nabla }_{\vec{t}}\hspace{2mm} v) =
\mu_{*}(-A_{v}\vec{t} + D^{\perp }_{\vec{t}}\hspace{2mm} v) = \mu_{*}(-A_{v}\vec{t})
$$

\noindent which is tangent to $(\mu \circ \gamma )$, thus proving $D^{\perp }_{\mu_{*}\vec{t}}\hspace{2mm}
\mu_{*}v=0$.

\bigskip

(2) It's a direct consequence of part (1) and claims 1 and 2 at the beginning of the proof.$\Box$

\bigskip

\textbf{Acknowledgments} The author wants to express his gratitude to his colleagues Marco Castrill{\'{o}}n,
Laureano Gonz{\'{a}}lez-Vega, Bert J\"{u}tler and Gema Quintana for the useful talks about the theory of RM
vectors and frames.

\end{document}